\documentclass[]{amsart}
\usepackage{amsmath,amsfonts,amssymb}
\usepackage{verbatim}
\usepackage{url}
\def\N{\mathbb N}

\def\Z{\mathbb Z}

\def\F{\mathbb F}
\def\ord{\mathop{\rm ord}}
\def\gcd{\mathop{\rm gcd}}

\theoremstyle{plain}
\newtheorem{theorem}{Theorem}[section]

\newtheorem{remark}[theorem]{Remark}
\newtheorem{example}[theorem]{Example}
\def\proof{{\it Proof: }}
\def\qed{\hfill\hbox{$\square$}}

\theoremstyle{definition}

\author[F.E. Brochero Mart\'{\i}nez]{F. E. Brochero Mart\'{\i}nez}
\author[C. R. Giraldo Vergara]{C. R. Giraldo Vergara}
\address{
Departamento de Matem\'{a}tica\\
Universidade Federal de Minas Gerais\\
UFMG\\
Belo Horizonte, MG\\
 30123-970\\
 Brazil\\
 }
 \email{fbrocher@mat.ufmg.br }\email{carmita@mat.ufmg.br}

\date{\today
}

\subjclass[2000]{}

\subjclass[2010]{16S34(primary) and 94B05(secondary)} 
\title{Explicit Idempotents of  finite group algebras}
\keywords{Irreducible cyclic codes, Primitive Idempotents}
\begin{document}

\begin{abstract}
Let $\F_q$ be a finite field, $G$  a finite cyclic group of order $p^k$ and $p$ is an odd prime with $\gcd(q,p)=1$.  
In this article, we determine an explicit expression for the primitive idempotents of $\F_qG$. 
This result extends the result in   \cite{ArPr1},  \cite{ArPr2}  and \cite{SBDR}. 
\end{abstract}

\maketitle

\section{Introduction}

Let $G$ be a finite cyclic group of order $n$ and $\F_q$  a finite field of order $q$, where $q$ is prime relative to $n$.  The
 cyclic codes of length $n$ over $\F_q$ can be viewed as an ideal in the group algebra $\F_qG$ and each
 ideal is generated by an idempotent of $\F_qG$.  By the representation theorem of abelian groups we know   that
$$G\simeq C_{p_1^{\beta_1}}\times \cdots \times C_{p_r^{\beta_r}}$$
where $C_{p_j^{\beta_j}}$ is a cyclic group of order $p_j^{\beta_j}$ and  $p_1,\dots, p_r$ are  distinct primes. In addition, it is well known  that
$$\F_qG\simeq  \F_qC_{p_1^{\beta_1}}\otimes \cdots \otimes \F_qC_{p_r^{\beta_r}}.$$
From this fact, in order to construct the idempotents of the  cyclic group algebra $\F_qG$, it is enough to consider  the case $G=C_n$ where  $n$ is a power of a prime.
Observe that the condition $\gcd(n,q)=1$ is necessary by the Maschke theorem (see \cite{CuRe} theorem 10.8).

\section{Primitive idempotents: General calculation}
Let $\Phi_d(x)$ denote the $d$-th cyclotomic polynomial, i.e.,  $\Phi_d(x)$ can be defined recursively by $\Phi_1(x)=x-1$ and $x^k-1=\prod_{d|k} \Phi_d(x)$.  It is well known (see \cite{LiNi} page 65 theorem 2.47) 
 that if $\gcd(q,d)=1$ then $\Phi_d(x)$ can be factorized into $r_d=\frac{\varphi(d)}{s_d}$ distinct monic irreducible polynomials of the same degree $s_d$ over $\F_q$ and $s_d=\ord_d q=\min\{k\in \N^*| q^k\equiv 1\pmod d\}$, i.e. 
$\Phi_d(x)$ can be factorized in $\F_q[x]$  as $f_{d,1}\cdot f_{d,2}\cdots f_{d,r_d}$,  where  each $f_{d,j}$ is an irreducible polynomial of degree $s_d$, and then 
$$x^n-1=\prod_{d|n} \prod_{j=1}^{r_d} f_{s,j}.$$
Observe that if $K$ is a decomposition field of the cyclotomic polynomial $\Phi_d(x)$,  then for each  pair $f_{d,i}, f_{d,j}$ there exists $\tau\in Gal(K|\F_q)$ such that $\tau(f_{d,i})= f_{d,j}$.

By the Chinese remainder theorem, we know that
$$\F_qC_n\simeq \frac {\F_q[x]}{\langle x^n-1\rangle}\simeq \bigoplus_{d|n}\bigoplus_{j=1}^{r_d}
\frac {\F_q[x]}{\langle f_{d,j}\rangle}
$$
where the  $\F_q$-algebra isomorphisms are  naturally defined  using a  generator $g$ of $C_n$ as  
$g\mapsto \overline x\mapsto (\overline x,\dots, \overline x)$.

Since each direct sum term is a field, then this decomposition is a Weddeburn decomposition of the group algebra and each primitive idempotent is of the form $(\overline 0,\dots,\overline 0,\overline 1,\overline 0, \dots, \overline 0)$.
Therefore, if  $ e_{d,j}$ is a primitive idempotent of $\F_qC_n$, then it can be seen as a polynomial $e_{d,j}(x)$ with the following properties :
\begin{enumerate}
\item  $deg(e_{d,j}(x))<n$
\item $e_{d,j}(x)$ is divisible by $f_{d_1,j_1}$  for all $(d_1,j_1)\ne (d,j)$
\item  $e_{d,j}(x)-1$ is divisible by $f_{d,j}$.
\end{enumerate}
From these three properties, we have the following
\begin{theorem}\label{faztudo} Let $\F_q$ be a finite field with $q$ element and $n\in \N^*$ such that $gcd(q,n)=1$, then  each primitive idempotent of $\F_qC_n$ is of the form
 $$e_{d,j}(x)=\frac {x^n-1}{f_{d,j}(x)} h_{d,j}(x),$$
where $f_{d,j}$ is an irreducible factor of the cyclotomic polynomial $\Phi_d(x)$, $d$ is a divisor of $n$
and $h_{d,j}\in \F_q[x]$ is a polynomial with $\deg(h_{d,j})<s_d:=\ord_d q$  that   is the inverse of $\frac {x^n-1}{f_{d,j}(x)}$ in the field $\frac {\F_q[x]}{\langle f_{d,j}\rangle}$. 
\end{theorem}

Observe that if we know the polynomial $f_{d,j}$ then   $h_{d,j}$ can be explicitly calculated using the Extended Euclidean Algorithm for polynomials.  In general,  the factorization of $\Phi_d(x)$ in $\F_q[x]$ for arbitrary $d$ and $q$ is an open problem. Some especial cases can be found in \cite{LiNi}, \cite{FiYu} and \cite{Mey}.

\section{The idempotents for some special known cases}
In this section we are going to show, without proof,  the idempotents in the extremal cases where $x^n-1$ factorized into linear factors in $\F_q[x]$ and where  each  cyclotomic polynomial $\Phi_d(x)$, with $d|n$, is irreducible in $\F_q[x]$.
 Before that, we need the following remarks:

\begin{remark}\label{obs1} The cyclotomic polynomial $\Phi_d(x)\in \F_q[x]$, with $\gcd(d,q)=1$ is factorized into linear factors if and only if $q\equiv 1 \pmod d$, thus, $x^n-1\in \F_q[x]$ is factorized into linear factors if and only if $q\equiv 1 \pmod n$.
\end{remark}

\begin{remark}\label{obs2} The cyclotomic  polynomial $\Phi_d(x)\in \F_q[x]$, with $\gcd(d,q)=1$  is irreducible if and only if  $q$ is a primitive root modulus $d$, i.e., $\ord_d q=\varphi(d)$.
\end{remark}

\begin{remark}\label{obs3} The group $\Z_d^*$ has a primitive root if and only if 
 $d=2,4,p^i$ or $2p^i$, where $p$ is an odd prime and $i\in \N^*$. 
\end{remark}

\begin{theorem}[\cite{ArPr1} theorem 2.1]\label{TAP1}
Let $\F_q$ be a finite field, and $n\in \N^*$ such that $q\equiv 1 \pmod n$. Then the primitive idempotents of $\dfrac {\F_q[x]}{\langle x^n-1\rangle}$ are given by
$$e_j(x)=\frac 1n \sum_{l=0}^{n-1} \zeta_n^{-jl} x^l, \quad 0\le j\le n-1,$$ 
where $\zeta_n\in \F_d$ is an $n$-th primitive  root of unity.
\end{theorem}

\begin{theorem}[\cite{ArPr1} theorem 3.5]\label{TAP2}
Let $\F_q$ be a finite field and $n=p^k$, where $p$ is an odd prime and $k\ge 1$ such that 
$\ord_{n} q=\varphi(n)$. Then the primitive idempotents of $\dfrac {\F_q[x]}{\langle x^n-1\rangle}$ are given by
$$e_0(x)=\frac 1{p^k}\sum_{l=0}^{p^k-1} x^l$$
and 
$$e_j(x)=\frac 1{p^{k-j+1}}\sum_{l=0}^{p^{k-j+1}-1} x^{p^{j-1}l }-\frac 1{p^{k-j}}\sum_{l=0}^{p^{k-j}-1} x^{p^jl} \quad 1\le j\le k.$$
\end{theorem}
Observe that the representation shown here  of the idempotents is the same one 

found by Ferraz and Polcino Milies in \cite{FePo}.

\section{The case $p|( q-1)$}

Let $\F_q$ be a finite field such that $p^{m}||(q-1)$ (i.e, $p^m|(q-1)$ and $p^{m+1}\nmid (q-1)$),  where $p$ is an odd prime and $m\ge 1$.  It follows that there exists   a primitive $p^{m}$-th root of unity $\zeta_{p^m}$ in $\F_q$.  In addition,
$s_j=\ord_{p^j} q=\begin{cases} 1 &\hbox{if } j\le m\\
p^{j-m}&\hbox{if } j>m\end{cases}$ and, therefore, $\Phi_{p^j}(x)$ is  factorized into linear factors if $j\le m$ and, in the case $j>m$ the factorization in irreducible factors is
$$\Phi_{p^j}(x)=\Phi_{p^m}(x^{p^{j-m}})=\prod_{{l=1}\atop{ (p,l)=1}}^{p^{m}-1} (x^{p^{j-m}}- \zeta_{p^m}^l).$$
We observe that in this factorization is essential that $p$ be an odd prime, or in the case $p=2$ it is necessary  that $q\equiv 1 \pmod 4$.
Using this fact, we obtain the following new result:

\begin{theorem}\label{case_q-1} Let $\F_q$ be a finite field and $p$  is a prime such that $p^m||\gcd(p^k,  q-1)$, where $m\ge 1$. If $q\equiv 1 \pmod 4$ or $p$ is an odd prime, then the primitive idempotents of $\dfrac {\F_q[x]}{\langle x^{p^k} -1\rangle }$ are of the following forms:
\begin{enumerate} 
\item  For each $0\le j\le p^m-1$, 
$$ \displaystyle e_j(x)= \frac 1{p^k} \sum_{l=0}^{p^k-1} \zeta_{p^m}^{-jl} x^l.$$
\item For each $m< s\le k$ and $0< l<p^m$ such that $\gcd(l,p)=1$
$$e_{s,l}(x)=\frac  {1}{p^{k+m-s}}\sum_{j=0}^{p^{k-s}-1} \zeta_{p^m}^{-lj} x^{p^{k+s-m}j}$$
\end{enumerate}

\end{theorem}
\proof Observe that 
$$\frac {\F_q[x]}{\langle x^{p^k}-1\rangle}\simeq \frac  {\F_q[x]}{\langle x^{p^m}-1\rangle}\oplus\bigoplus_{s=m+1}^k  \frac  {\F_q[x]}{\langle \Phi_{p^s}(x)\rangle},$$
where the second summand  does not appear in the case $k=m$.  
Therefore, by  Remark \ref{obs1}, the first case corresponds  to the idempotents  associated to the factor $x-\zeta_{p^m}^j\in \F_q[x]$ of $x^{p^k}-1$, and  by  theorem  \ref{faztudo} we know that the
idempotent associated to this factor is of the form
$$e_j(x)=P_j(x) h_j,$$
where $P_j(x)=\frac {x^{p^m}-1}{x- \zeta_{p^m}^j}$ and  $h_j\in \F_q$ are such that $P_j(\zeta_{p_m}^j)h_j=1$. Observe that
$$P_j(x)=\frac {x^{p^m}-1}{x- \zeta_{p^m}^j}=\frac {x^{p^m}-(\zeta_{p^m}^j)^{p^m}}{x- \zeta_{p^m}^j}=\sum_{l=0}^{p^m-1} x^{l}\zeta_{p^m}^{j(p^m-1-l)}= \sum_{l=0}^{p^m-1} \zeta_n^{-jl-j} x^l,$$
and 
 $P_j(\zeta_{p^m}^j)=
p^m\zeta_{p^m}^{-j}$. Thus $h_j=\frac {\zeta_{p^m}^{j}}{p^m}$ and this implies the result of the first case of the theorem.

For the second case, let $x^{p^{s-m}}- \zeta_{p^m}^l$ be an irreducible factor of $\Phi_{p^s}(x)$, then the associated primitive idempotents  are
$$e_{s,l}= P_{s,l}(x)h_{s,l}(x),$$
where $P_{s,l}(x)=\dfrac {x^{p^k}-1}{x^{p^{s-m}}- \zeta_{p^m}^l}$,  and $h_{s,l}(x)$ is a polynomial satisfying
$$P_{s,l}(x)h_{s,l}(x)\equiv 1 \pmod {x^{p^{s-m}}- \zeta_{p^m}^l}\quad\hbox{ and }\quad\deg(h_{s,l}(x))< p^{s-m}.$$
Substituting $x^{p^{s-m}}$ by  $y$, it follows that $P_{s,l}(x)=\tilde P(y)=\dfrac {y^{p^{k+m-s}}-1}{y-\zeta_{p^m}^l}$, and
using a formal version of L'H\^opital rule, we obtain
$$\tilde P(y)\equiv P(1)=p^{k+m-s} \zeta_{p^m}^{l(p^{k+m-s}-1)}=p^{k+m-s}\zeta_{p^m}^{-l} \pmod {y-\zeta_{p^m}^l}$$
or, equivalently,
$ P_{s,l}(x)\equiv p^{k+m-s}\zeta_{p^m}^{-l}\pmod {x^{p^{s-m}}-\zeta_{p^m}^l}$.

Then 
$h_{s,l}=\frac {\zeta_{p^m}^l} {p^{k+m-s}}$ and
\begin{align*}
e_{s,l}(x)&=\frac  {\zeta_{p^m}^l}{p^{k+m-s}}\dfrac {x^{p^k}-1}{x^{p^{s-m}}- \zeta_{p^m}^l}=\frac  {\zeta_{p^m}^l}{p^{k+m-s}}\dfrac {x^{p^k}-(\zeta_{p^m}^l)^{p^k}}{x^{p^{s-m}}- \zeta_{p^m}^l}\\
&=\frac  {1}{p^{k+m-s}}\sum_{j=0}^{p^{k-s+m}-1} \zeta_{p^m}^{l(p^{k-s+m}-j)} x^{jp^{s-m}}=\frac  {1}{p^{k+m-s}}\sum_{j=0}^{p^{k-s+m}-1} \zeta_{p^m}^{-lj} x^{jp^{s-m}}
\end{align*}
concluding the proof.
\qed

\begin{remark} The case when  $p^m||(q-1)$ with $m\ge k$, is a particular  case of theorem  \ref{TAP1}  when $n=p^k$.
\end{remark}

\section{General Case }
Let $\F$ be a finite field with $q$ element and $G$ a group $p^k$ element, where $\gcd(q,p)=1$. 
The classical method to calculate the irreducible idempotents depends of the computation of the irreducible characters $\psi:G\to \widehat \F$, where $\widehat \F$ denotes the algebraic closure of $\F$,   and the Galois group $Gal(\F(\psi),\F)$. In fact,  $e(\psi)=\frac1{p^k} \sum_{g\in G} \psi(g^{-1}) g$ is a primitive idempotent of 
$\widehat \F G$ and $$e_{\F}(\psi)=\sum_{\sigma \in Gal(\F(\psi),\F)} \sigma \cdot e(\psi)$$
is a primitive idempotent of $\F G$, where $\sigma$ acts on the coefficient of $e(\psi)$.

In this section, we are going to calculate the idempotent, without calculate the irreducible characters, only using the trace of some extension of $\F$. 

Suppose that  $\ord_p q=t>1$ and   $m$ is an integer such that $p^m||(q^t-1)$. By little Fermat theorem, it's known that $t|(p-1)$.  
Under such condition, $\F_{q^l} $ does not have a $p$-th primitive root of unit for all $l<t$, but 
$\F_{q^t}$ contains $\zeta_{p^m}$,  a primitive $p^m$-th root of unit,  then $\F_{q^t}$ can be seen as a decomposition field of the minimal polynomial of 
$\zeta_{p^m}$ under $\F_q$, i.e., there exists an irreducible polynomial $Q(x)\in \F_q[x]$ of degree $t$, such that $Q(\zeta_{p^m})=0$. 

In addition, if 
$\tau\in Gal(\F_{q^t},\F_q)$ is the Fr\"obenius automorphism $a\mathop{\longmapsto}\limits^{\tau} a^q$,
 then  $$\{\tau^{j}(\zeta_{p^m})| j=0,1,\dots, t-1\}=\{\zeta_{p^m}^{q^j}| j=0,1,\dots, t-1\}$$ is the set 
of conjugates of  $\zeta_{p^m}$ over $\F_q$. In general, for all  $a\in \F_{q^t}$,  the function 
$$\begin{array}{rccl}\sigma_1:& \F_{q^t}&\to& \F_q\\
&a&\mapsto&a+\tau(a)+\tau^{2}(a)+\cdots+\tau^{(t-1)}(a),
\end{array}$$ 
is well defined.

By the previous section, we show the explicit form of the primitive idempotents of $\dfrac{\F_{q^t}[x]}{\langle x^{p^k}-1\rangle}$. The next theorem uses this representation in order to calculate the form
of the primitive idempotents of $\dfrac{\F_q[x]}{\langle x^{p^k}-1\rangle}$.

\begin{theorem}\label{geral} Let $\F_q$ be a finite field and assume  $p$ is an odd prime such that $\ord_p q=t>1$  and  $p^m|| (q^t-1)$, where $ m\ge 1$. Then the primitive idempotents of $\dfrac {\F_q[x]}{x^{p^k} -1}$ are of the following forms:
\begin{enumerate} 
\item $ \displaystyle e_0(x)= \frac 1{p^k} \sum_{l=0}^{p^k-1} x^l.$
\item  For each $0< j\le p^m-1,$ 
$$ \displaystyle e_j(x)= \frac 1{p^k} \sum_{l=0}^{p^k-1} \sigma_1( \zeta_{p^m}^{-jl}) x^l,$$
where $e_{i}$ and $e_{j}$ are the same idempotents if and only if $i\equiv jq^u\pmod{p^m}$ for some $u\in\Z$, and,
therefore, there are $\frac {p^m-1}{t}$  different primitive idempotents of this type.

\item For each $m< s\le k$ and $0< l<p^m$ such that $\gcd(l,p)=1$,
$$e_{s,l}(x)=\frac  {1}{p^{k+m-s}}\sum_{j=0}^{p^{k-s+m}-1} \sigma_1\bigl(\zeta_{p^m}^{-lj)}\bigr )x^{p^{s-m}j},$$
where $e_{s,l_1}$ and $e_{s,l_2}$ are the same idempotents if and only if $l_1\equiv l_2q^u\pmod{p^m}$, for some $u\in\N$, and, therefore, for each $s$ fixed, there are $\frac {\varphi(p^{k-s+m})}{t} $ different primitive idempotents of this type.
\end{enumerate}

\end{theorem}
\begin{remark}  In \cite{SBDR}, using a cyclotomic  cosets method,  was studied the particular case when  $m=1$.
\end{remark}
\proof
Let $E(x)$ be a primitive idempotent of $\dfrac{\F_q[x]}{\langle x^{p^k}-1\rangle}$. It follows that $E(x)$ is also an idempotent of $\dfrac{\F_{q^t}[x]}{\langle x^{p^k}-1\rangle}$ and therefore $E(x)$ is a direct sum of primitive idempotents of $\dfrac{\F_{q^t}[x]}{\langle x^{p^k}-1\rangle}$.  In the case that $E(x)$ be also primitive in $\dfrac{\F_{q^t}[x]}{\langle x^{p^k}-1\rangle}$, then by theorem   \ref{case_q-1} we know  the unique idempotent  with this propriety is  $e_0(x)$, i.e., the case  when  $j=0$ and therefore $E(x)=e_0(x)$.

Now, suppose that $E(x)\ne e_0(x)$, and let $e(x)\in \dfrac{\F_{q^t}[x]}{\langle x^{p^k}-1\rangle}$ be a primitive idempotent such that $e(x)\cdot E(x)=e(x)\notin\dfrac{\F_q[x]}{\langle x^{p^k}-1\rangle} $.  Since $\tau (E(x))=E(x)$, it follows that $\tau(e(x))$ is also a direct summand of $E(x)$. In addition,  it is known that
 $\tau^{r}(\zeta_{p^m}^j)=\zeta_{p^m}^j$  if and only if 
$j(q^r-1)$  is divisible by $p^m$ and this is equivalent  to $p^m|j$ or $t|r$.   From this, we conclude that 
$\tau^{r}(e(x))=e(x)$ if and only if  $t|r$, thus   
$$\{e(x),\tau(e(x)),\dots, \tau^{(t-1)}(e(x))\}$$
is a list of different idempotents that are direct summands of $E(x)$. Finally, since $\sigma_1(e(x)):=e(x)+\tau(e(x))+\cdots +\tau^{(t-1)}(e(x))\in \dfrac{\F_q[x]}{\langle x^{p^k}-1\rangle}$ and using the fact that $E(x)$ is primitive,  we conclude that 
$E(x)=\sigma_1(e(x))$,  in other words, we can obtain every primitive idempotent of $\dfrac{\F_q[x]}{\langle x^{p^k}-1\rangle}$ different of $e_0(x)$ from the idempotent of $\dfrac{\F_{q^t}[x]}{\langle x^{p^k}-1\rangle}$ and the ring homomorphism $\sigma_1$. 
Thus, the other cases of the theorem follow directly from the cases (1) and (2) of  Theorem \ref{case_q-1}.\qed

\section{Sage implementation and examples}
In this section, some examples are shown explicitly. In order to  find these idempotents, we have implemented the last theorem  in the SAGE program\footnote{\url{http://www.sagemath.org}}, as it is shown in the following code:

First we defined the field $\F_{q^t}$, the $p^m$-th root of the unity, and the polynomial ring $\F_{q^t}[x]$
\begin{verbatim}
sage: k.<a>=GF(q^t,'a'); 
sage: b=a^((q^t-1)/p^m); 
sage: F.<x>=PolynomialRing(k,'x')
\end{verbatim}

Implementation the function $\sigma_1(\zeta_{p^m}^{lu})$ 
\begin{verbatim}
sage: def sigma(l,u):
...       sumconj=sum([b^(l*u*q^i) for i in range(0,t)]);
...       return(sumconj)
\end{verbatim}

The idempotents of the second type
\begin{verbatim}
sage: def Idemp2(l):
...       v = []
...       for i in range(0,p^k):
...          v.append( sigma(q,t,l,i)/p^k)
...       Poli= sum([v[j]*x^j for j in range(0,p^k)])
...       return Poli
\end{verbatim}

The idempotents of the third type
\begin{verbatim}
sage: def Idemp3(l,s):
...       v = []
...       for j in range(0,p^(k-s+m)):
...          v.append( sigma(q,t,l,p^(k-s+m)-j)/p^(k+m-s))
...       Poli= sum([v[j]*x^(p^(s-m)*j) for j in range(0,p^(k-s+m))])
...       return Poli
\end{verbatim}

\begin{example} 
In the group ring $\F_{17}C_{13^2}$, since $13||17^6-1$, then  $t=6$ and $m=1$. Thus,  there exists $\frac {13-1}{6}=2$ primitive idempotents of the second type, and making $s=2=k$, there exists $\frac {\varphi(13)}6=2$    primitive idempotents of the third type. In fact, the primitive idempotents are:
\begin{itemize}
\item $\displaystyle e_0=16 \sum\limits_{i=0}^{168} x^i
$;
\item two idempotents  of the second type
\begin{align*}
e_{1}=\sum\limits_{i=0}^{12} x^{13i}&(5 x^{12} + 13
x^{11} + 5 x^{10} + 5 x^{9} + 13 x^{8} + 13 x^{7} + 13 x^{6} +\\
&+ 13 x^{5}
+ 5 x^{4} + 5 x^{3} + 13 x^{2} + 5 x + 11)\\
&\\
e_2=\sum\limits_{i=0}^{12} x^{13i}&(13 x^{12} + 5x^{11} + 13 x^{10} + 13 x^{9} + 5 x^{8} + 
5 x^{7} + 5 x^{6} +\\&+ 5 x^{5} +
13 x^{4} + 13 x^{3} + 5 x^{2} + 13 x + 11)
\end{align*}
\item 
and two idempotents of the third type
\begin{align*}
e_{2,1}=&        	
14 x^{156} + 16 x^{143} + 14 x^{130} + 14 x^{117} + 16 x^{104}
+ 16
x^{91} + 16 x^{78}+\\ &+ 16 x^{65} + 14 x^{52} + 14 x^{39} + 16 x^{26} + 14
x^{13} + 7\\
&\\
e_{2,2}=&16 x^{156} + 14 x^{143} + 16 x^{130} + 16 x^{117} + 14 x^{104} + 
14
x^{91} + 14 x^{78} +\\ &+ 14 x^{65} + 16 x^{52} + 16 x^{39} + 14 x^{26} + 16
x^{13} + 7
\end{align*}
\end{itemize}
\end{example}

\end{document}